\documentclass[a4paper,11pt]{article}
\usepackage{xypic,fancyhdr}
\usepackage{geometry,citesort,enumerate,amsmath,amssymb,latexsym,theorem}
 \xyoption{curve}

{\theorembodyfont{\rmfamily}\newtheorem{example}{Example}[section]}
{\theorembodyfont{\rmfamily}}
\newtheorem{prop}[example]{Proposition}
\newtheorem{thm}[example]{Theorem}
{\theorembodyfont{\rmfamily}}
\newtheorem{cor}[example]{Corollary}

\def\<{\langle}
\def\>{\rangle}
\def\Z{\mathbb{Z}}
\def\lan{\langle}
\def\ran{\rangle}

\newcommand{\sast}{_{\ast}}

\newcommand{\gpd}{\mathsf{Gpd}}

\def\crs{\mathsf{Crs}}
\def\Crs{\mathsf{CRS}}
\def\ftop{\mathsf{FTop}}

\def\cM{\mathcal{M}}

\newcommand{\Env}[2]{\begin{#1} #2\end{#1}}

\newcommand{\eps}{\varepsilon}
\newcommand{\io}{^{-1}}

\def\ge{\geqslant}

\def\ges{\geqslant}
\def\leq{\leqslant}

\newcommand{\labto}[1]{\stackrel{#1}{\longrightarrow}}

\def\phi{\varphi}

\def\E{\mathsf{E}}
\def\sC{\mathsf{C}}
\newcommand{\cP}{{\mathcal P}}
\newcommand{\I}{{\mathcal I}}
\def\subs{\subseteq}
\newenvironment{proof}{\noindent {\bf Proof} }{ \hfill
$\Box$ \mbox{}}

\def\gpd{\mathsf{Gpd}}
\def\Gpd{\mathsf{GPD}}

\def\cM{\mathcal{M}}

\def\epsilon{\varepsilon}
\begin{document}\title{Crossed complexes, free crossed resolutions and  \\
graph products of groups}
\author{R. Brown \thanks{\{ R. Brown, T.Porter\}, Mathematics Division,  School of Informatics,
University of Wales, Bangor,Gwynedd LL57 1UT,
 U.K. email: \{r.brown, t.porter\}@bangor.ac.uk} \and M.
 Bullejos
 \thanks{Departmento da Algebra, Universidad, Granada, Spain. email: bullejos@ugr.es} \and T.
Porter$\,^*$}

 \maketitle
\thispagestyle{plain}

\begin{abstract} The category of crossed
complexes gives an algebraic model of $CW$-complexes and cellular
maps. Free crossed resolutions of groups contain information on a
presentation of the group as well as higher homological
information. We relate this to the problem of calculating
non-abelian extensions. We show how the strong properties of this
category allow for the computation of free crossed resolutions of
graph products of groups, and so obtain computations of higher
homotopical syzygies in this case.  \footnote{KEYWORDS: crossed
complex,
resolution, higher syzygies, graph product. \\
\rule{1.5em}{0mm}MATHSCI 2000 CLASS: 55U99, 20J99, 18G55}
\end{abstract}
\section*{Introduction}

One  aim of this paper is to advertise the category of crossed
complexes, and the notion of  free crossed resolution, as a
working tool for certain problems in combinatorial group theory.
Crossed complexes are an analogue of chain complexes of modules
over a group ring, but with a non abelian part,  a {\it crossed
module},  at the bottom dimensional part. This allows for crossed
complexes to contain in that part the data for a presentation of a
group, and to contain in other parts  higher homological data. The
non abelian nature, and also the generalisation to groupoids
rather than just groups, allows for a closer representation of
geometry, and this, combined with very convenient properties of
the category of crossed complexes, can allow for more and easier
calculations than are available in the standard theory of chain
complexes of modules.

The notion of {\it crossed complex} of groups was defined by
A.L.Blakers in 1946 \cite{Blakers} (under the term `group system')
and Whitehead \cite{W1}, under the term `homotopy system' (except
that he restricted to the free case). Blakers used these as a way
of systematising known properties of relative homotopy groups
$\pi_n(X_n,X_{n-1},x)$ of a filtered space
$$X_*: X_0 \subs X_1 \subs X_2\subs \cdots \subs X_n \subs \cdots
\subs X_{\infty} .$$ It is significant that he used the notion to
establish relations between homotopy and homology of a space.
Whitehead was strongly concerned with realisability, that is with
the passage between algebra and geometry and back again. He
explored the relations between crossed complexes and chain
complexes with a group of operators and established  remarkable
realisability properties, some of which we explain later.

There was another stream of interest in crossed complexes, but  in
a broader algebraic framework, in work of Frohlich \cite{Fro} and
Lue \cite{lue1}. This gave a general formulation of cohomology
groups relative to a variety in terms of equivalence classes of
certain exact sequences. However the relation of these
equivalence classes with the usual cohomology of groups was not
made explicit till papers of Holt \cite{Ho} and Huebschmann \cite
{Hue}. The situation is described in Lue's paper \cite{lue2}.

Since our interest is in the relation with homotopy theory, we are
interested in the case of groups rather than other algebraic
systems. However there is one key change we have to make, as
stated above, namely that we have to generalise to groupoids
rather than groups. This makes for a more effective modelling of
the geometry, since we need to use $CW$-complexes which are non
reduced, i.e. have more than one $0$-cell, for example universal
covering spaces, and simplices. It allows for easier statement of
some theorems because the coproduct for groupoids is simply
disjoint union. It also gives the category of crossed complexes
better algebraic properties, principally that it is a monoidal
closed category in the sense of having an internal hom which is
adjoint to a tensor product. This is a generalisation of a
standard property of groupoids: if $\gpd$ denotes the category of
groupoids, then for any groupoids $A,B,C$ there is a natural
bijection, the exponential law,
$$\gpd(A \times B,C) \cong \gpd(A,\Gpd(B,C))$$ where $A \times B$
is the usual product of groupoids, and $\Gpd(B,C)$ is the groupoid
whose objects are the morphisms $B \to C$ and whose arrows are the
natural equivalences (or conjugacies) of morphisms. Another
advantage of groupoids is that there is a standard model of the
unit interval in topology, namely the groupoid $\I$ with two
objects $0,1$ and exactly one  arrow $\iota: 0 \to 1$. This leads
to a homotopy theory for groupoids in terms of morphisms $\I
\times B \to C$ (in the case of groups, homotopies of morphisms
are just conjugacies). It also leads to a useful notion of {\em
fibration} of groupoids (see for example \cite{rbbook}).

This exponential law for groupoids  is modelled in the category
$\crs$ of crossed complexes by a natural isomorphism $$ \crs(A
\otimes B,C) \cong \crs(A, \Crs(B,C));$$ that is, $\crs$ is a
monoidal closed category, as proved by Brown and Higgins in
\cite{BHtens}. The groupoid $\I$ determines a crossed complex also
written $\I$ and so a homotopy theory for crossed complexes in
terms of a cylinder object $\I \otimes B$ and homotopies of the
form $\I \otimes B \to C$. For our purposes, the key result is the
tensor product $A \otimes B$; this has a complicated formal
definition, reflecting the algebraic complexity of the definition
of crossed complex. However, for the purposes of calculating with
free crossed complexes, it is sufficient to know the boundaries of
elements of the free bases, and also the value of morphisms on
these elements. Thus the great advantage is that the free crossed
resolutions model Eilenberg-Mac~Lane spaces and their cellular
maps (see Corollary \ref{cwmodel} and Proposition \ref{cwmaps}),
and give modes of calculating with these which would be very
difficult geometrically.

The end product of this paper (section 4) is to show how these
methods enable one to compute higher homotopical syzygies for
graph products of groups, generalising to higher dimensions the
results of \cite{BHP}. A further paper \cite{EMTW} will apply
related methods to obtain free crossed resolutions for amalgamated
sums and HNN-extensions of groups.

\section{Definitions and basic properties}\label{algebra}
A {\em crossed complex } $C$ (of groupoids) is a sequence of
morphisms of groupoids over $C_0$ $$\diagram \cdots \rto & C_n
\dto<-.05ex>^{\beta}\rto^-{\delta_n} & C_{n-1} \rto
\dto<-1.2ex>^{\beta} & \cdots \rto &
 C_2\rto^{\delta_2} \dto<-.05ex>^{\beta} & C_1
\dto<0.0ex>^(0.45){\delta^1} \dto<-1ex>_(0.45){\delta^0} \\ &
C_0&C_0\rule{0.5em}{0ex}  & & \rule{0.5em}{0ex} C_0 &
\rule{0em}{0ex} C_0 . \enddiagram$$ Here $\{ C_n\}_{n \ge 2}$ is
a family of groups with base point map $\beta$, so that  for $p
\in C_0, $ we have groups $C_n(p)= \beta \io (p)$,
 and $ \delta^0, \delta^1$ are the source and targets
for the groupoid $C_1$. We further require given an operation of
the groupoid $C_1$ on each family of groups $C_n$ for $n
\geqslant 2$ such that:
 \begin{enumerate}[(i)]
      \item  each $\delta_n$ is a morphism over the identity on $C_0$;

      \item  $C_2 \rightarrow C_1$ is a crossed module over $C_1$;

      \item  \label{modul} $C_n$ is a $C_1$-module for $n\ges 3$;

      \item  $\delta : C_n \rightarrow C_{n-1}$ is an operator morphism
 for $n\ges 3$;

      \item  $\delta\delta :C_n \rightarrow C_{n-2}$ is trivial for
$n\ges 3$;

      \item  $\delta C_2$ acts trivially on $C_n$ for $n\ges 3$.
   \end{enumerate}
Because of axiom (iii) we shall write the composition in $C_n$
additively for $n \ges 3$, but we will use multiplicative notation
in dimensions 1 and 2 (except when giving the rules for the
tensor product). Note that if $a:p\to q, b:q \to r$ in $C_1$ then
the composition is written $ab: p \to r$. If further $x \in
C_n(p)$ then $x^a \in C_n(q)$ and the usual laws of an action
apply. We write $C_1(p) = C_1(p,p)$, and $C_1$ operates on this
family of groups by conjugation.  The condition (iv) then implies
that $\delta_2(x^a)= a \io \delta_2(x) a$, while condition (ii)
gives further that $x\io y x= y^{\delta_2(x)}$ for $x,y \in
C_2(p), a \in C_1(p,q)$. Consequently $\delta_2(C_2) $ is normal
in $C_1$, and $Ker \;\delta_2 $ is central in $C_2$ and is
operated on trivially by $\delta_2(C_2) $.

Let $C$ be a crossed complex. Its {\em  fundamental groupoid}
$\pi_1{C}$ is the quotient of the groupoid $C_1$ by the normal,
totally disconnected subgroupoid $\delta{C_2}$. The rules for a
crossed complex  give $C_n$, for $n\ges 3$, and also $Ker
\;\delta_2 $,  the induced structure of $\pi_1{C}$-module.

The crossed complex $C$ is {\em reduced } if $C_0$ is a singleton,
so that all the groupoids $C_n, n \ge 1$ are groups. This was the
case  considered in \cite{Blakers,W1} and many other sources.

A {\it morphism} $f: C \rightarrow D$ of crossed complexes is a
family  of groupoid morphisms $f_n : C_n \rightarrow D_n ~(n\ges
0)$ which preserves all the structure. This defines the category
${\crs}$ of crossed complexes. The fundamental groupoid now gives
a functor $ \pi_1 : \crs \to \gpd $. This functor is left adjoint
to the functor $i : \gpd \to  \crs $ where for a groupoid $G$ the
crossed complex  $iG$ agrees with $G$ in dimensions 0 and 1, and
is otherwise trivial.

An $m$-{\em truncated} crossed complex $C$ consists of all the
structure defined above but only for $n \leq m$. In particular, an
$m$-truncated crossed complex is for $m=0,1,2$ simply a set, a
groupoid, and a crossed module respectively.

One basic algebraic example of crossed complex comes from the
notion of {\it identities among relations} (For more details on
the following, see \cite{RB&JH}.) Let $\cP=\lan X_1|w\ran$ be a
presentation of a group $G$ where $w$ is a function from a set
$X_2$ to $F(X_1)$, the free group on the set $X_1$ of generators
of $G$. This gives an epimorphism $\phi:F(X_1)  \to G$, with
kernel $N(R)$, the normal closure in $F(X_1)$ of the set
$R=w(X_2)$.

Let $H$ be the free $F(X_1)$-operator group on the set $X_2$, so
that $H$ is the free group on the elements  $(x,u) \in X_2 \times
F(X_1)$. Let $\delta': H \to F(X_1)$ be determined by $(x,u)
\mapsto u \io (wx) u$, so that the image of $\delta'$ is exactly
$N(R)$. Note that $F(X_1)$ operates on $H$ by $(x,u)^v= (x,uv)$,
and for all $h \in H, u \in F(X_1)$ we have \\
CM1) $\delta'(h^u)= u \io \delta'(h) u.$ \\
We say that $\delta': H \to F(X_1)$ is a {\it precrossed module}.

We now define {\it Peiffer commutators} for $h,k \in H$
$$\lan h,k\ran= h\io k\io h k^{\delta' h} .$$
Then $\delta'$ vanishes on Peiffer commutators. Also the subgroup
$P=\lan H,H\ran$ generated by the Peiffer commutators is a normal
$F(X_1)$-invariant subgroup of $H$. So we can define $C(w)= H/P$
and obtain the exact sequence
$$ C(w) \labto{\delta_2} F(X_1) \labto{\phi} G\to 1.$$

The morphism $\delta_2$ satisfies \\
CM2) $ c\io d c = d^{\delta_2 c}$ for all $c,d \in C(w)$. \\
The rules CM1, CM2) are the laws for a {\it crossed module}, and
$\delta_2 :C(w) \to F(X_1)$ is known as the {\it free crossed
$F(X_1)$-module on $w$}. The injection $i: X_2 \to C(w)$ has the
universal property that if $\mu :M \to F(X_1)$ is a crossed
module and $v: X_2 \to M$ is a function such that $\mu v =w$,
then there  is a unique crossed module morphism $\eta: C(w) \to
M$ such that $\eta i = w$. The elements of $C(w)$ are `formal
consequences'
$$c=\prod_{i=1}^n (x_i^{\eps_i})^{u_i}$$
where $n \ge 0,x_i \in X_2, \epsilon_i=\pm 1, u_i \in F(X_1)$,
$\delta_2(x^{\eps})^u=u\io (w x)^{\eps} u$, subject to the crossed
module rule $cd=dc^{\delta_2d}, c,d \in C(w)$.

The kernel $\pi(\cP)$ of $\delta_2$ is abelian and in fact
obtains the structure of $G$-module -- it is known as the
$G$-module of {\it identities among relations} for the
presentation.

We can now splice to the free crossed module any resolution of
$\pi(\cP)$ by free $G$-modules, and so obtain what is called a
{\it free crossed resolution} of the group $G$.

This construction is analogous to the usual construction of
higher order syzygies and free resolutions for modules, but taking
into account the non abelian nature of the group and its
presentation, and in particular the action of $F(X_1)$ on $N(R)$.

There is a notion of homotopy for morphisms of crossed complexes
defined using the tensor product and the crossed complex $\I$.
Assuming this we can state one of the basic homological results,
namely the uniqueness up to homotopy equivalence of free crossed
resolutions of a group $G$.

There is a {\it standard free crossed resolution} $F^{st}_*(G)$ of
a group $G$ in which $F^{st}_1(G)$ is the free group on the set
$G$ with generators $[a], a \in G$; $F^{st}_2(G)$ is the free
crossed $F^{st}_1(G)$-module on $w: G \times G \to F^{st}_1(G)$
given by
$$w(a,b)= [a][b][ab]\io , a,b \in G;$$
for $n \ge 3$, $F^{st}_n(G)$ is the free $G$-module on $G^n$, with
$$\delta_3[a,b,c]= [a,bc][ab,c]\io[a,b]\io [b,c]^{[a]\io},$$
and for $n \ge 4$
\begin{multline*}
\delta_n[a_1, a_2, \ldots, a_n]=[a_2, \ldots,a_n]^{a_1\io} +
 \sum _{i=1}^{n-1}(-1)^{i}
[a_1,a_2,\ldots,a_{i-1}, a_ia_{i+1},a_{i+2}, \ldots, a_n] +\\
\qquad \qquad +(-1)^{n} [a_1,a_2, \ldots, a_{n-1}].
\end{multline*}

There is an exact sequence
$$ \diagram F^{st}_2(G) \rto
^{\delta_2}  & F^{st}_1(G)\rto ^-{\phi} & G
\enddiagram$$ where $\phi([a])=a, \; a \in G$.

We can now see the advantage of this setup in considering the
notion of non abelian 2-cocycle on the group $G$ with values in
the group $K$. According to standard definitions, this is a pair
of functions $k^1: G \to Aut(K), k^2: G \times G \to K$
satisfying certain properties. But suppose $G$ is infinite; then
it is difficult to know how to specify these functions and check
the required properties.

However  the 2-cocycle definition turns out to be equivalent to
regarding $k^1,k^2$ as specifying a morphism of crossed complexes
$$ \diagram \cdots \rto & F^{st}_3(G) \rto ^{\delta_3} \dto &
F^{st}_2(G) \rto
^{\delta_2} \dto _{k^2} & F^{st}_1(G) \dto _{k^1} \\
\cdots \rto & 0 \rto & K \rto_-{\partial} & Aut(K)  \enddiagram
$$ (so   that $\partial k^2= k^1 \delta_2, \,  k^2 \delta_3=0$), where $
K \labto{\partial}  Aut(K)$ is the inner automorphism map and so
is a crossed module.  Further,  equivalent cocycles are just
homotopic morphisms.

Equivalent data to the above is thus obtained by replacing the
standard free crossed resolution by any homotopy equivalent free
crossed resolution.
\begin{example}
We let $\sC$ denote the infinite cyclic group. Then $\sC$ has a
free crossed resolution $1 \to  \cdots \to  \sC \labto{\phi} \sC$
with one free generator in dimension 1. It is shown in \cite{BW}
that the finite cyclic group of order $p$ with generator $t$, say,
$\sC_p$ has a free crossed resolution
$$ F(\sC_P): \cdots \to \Z[C_p] \labto{\delta_4} \Z[C_p] \labto{\delta_3} \Z[C_p]\labto{\delta_2}
C_\infty \labto{\phi} C_p $$ with a free generator $x_n$ in
dimension $n$ and $$\delta_n (x_n) = \begin{cases}
x_1^p &\text{if } n=2; \\
x_{n-1}(1-t)&\text{if } n \text{ is odd}; \\
x_{n-1}(1+t+t^2 + \cdots +t^{p-1})&\text{otherwise}.
\end{cases}$$ and where $\phi(x_1)=t$. It is shown in \cite{BP} how this enables one to
recover directly results of \cite{Z} on the enumeration and
classification of extensions by $\sC_p$. In fact a 2-cocycle on
$\sC_p$ can be specified (up to homotopy) by elements $k \in K,
a\in Aut(K) $ such that $\partial ( k) =a^n$ and $k=k^a$. The
extension $E$ determined by $k,a$ is given by $$E = (\sC \ltimes
K)/\{(t^n,k\io)\}$$   where the operation yielding the semidirect
product is given by $t$ operating via $a$.
\end{example}

Similar methods can be used to determine the 3-dimensional
obstruction class $l^3 \in H^3(G, A)$ corresponding to a crossed
module $\mu : M \to P$ with $Coker \; \mu = G, Ker \; \mu = A$,
provided we have a small free crossed resolution of the group $G$,
and so this focuses on methods of constructing such resolutions.
This method is successfully applied to the case with $G$ finite
cyclic in \cite{BW,BW2}.
\section{Relation with topology}\label{topology}
In order to give the  basic geometric example of a crossed complex
we first define a {\em filtered space} $X\sast $. By this we mean
a topological space $X_{\infty}$ and an increasing sequence of
subspaces $$  X \sast : X_0 \subseteq X_1 \subseteq  \cdots
\subseteq X_n \subseteq \cdots  \subseteq X_{\infty}.$$ A {\em
map} $f:X \sast \to Y\sast $ of filtered spaces consists of a map
$f:X_{\infty} \to Y_{\infty}$ of spaces such that for all $i \ges
0, f( X_i) \subseteq Y_i. $ This defines the category $\ftop$ of
filtered spaces and their maps. This category has a monoidal
structure in which $$(X_* \otimes Y_*)_n = \bigcup _{p+q=n}X_p
\times Y_q,$$ where it is best for later purposes to take the
product in the convenient category of compactly generated spaces,
so that if $X_*,Y_*$ are $CW$-spaces, then so also is $X_*
\otimes Y_*$.

We now define the {\em fundamental, or homotopy,  crossed complex
} functor
$$\pi : \ftop \to \crs.$$ If $C=\pi (X \sast) $, then $C_0=X_0,
$ and $C_1$ is the fundamental groupoid $\pi_1(X_1,X_0)$. For
$n\ges 2$, $C_n=\pi_n X \sast $ is the family of relative
homotopy groups $\pi_n(X_n,X_{n-1},p)$ for all $p \in X_0$. These
come equipped with the standard operations of $\pi_1 X \sast $ on
$\pi_n X \sast $ and boundary maps $\delta : \pi_n X \sast  \to
\pi_{n-1}X \sast $, namely the boundary of the homotopy exact
sequence of the triple $(X_n, X_{n-1},X_{n-2}).$ The axioms for
crossed complexes are in fact those universally satisfied by this
example, though this cannot be proved at this stage (see
\cite{Bhcolimits}).

This construction also explains why we want to consider crossed
complexes of groupoids rather than just groups. The reason is
exactly analogous to the reason for considering non reduced
$CW$-complexes, namely that we wish to consider covering spaces,
which automatically have more than one vertex in the non trivial
case. Similarly, we wish to consider covering morphisms of
crossed complexes as a tool for analysing presentations of
groups, analogously to the way covering morphisms of groupoids
were used for group theory applications by P.J. Higgins in 1964
in \cite{Hi1}.  A key tool in this is the use of paths in a Cayley
graph as giving elements of the free groupoid on the Cayley
graph, so that one moves to consider presentations of groupoids.
Further, as is shown by Brown and Razak in \cite{BR}, higher
dimensional information is obtained by regarding the free
generators of the universal cover of a free crossed resolution as
giving a higher order Cayley graph, i.e. a Cayley graph with
higher order syzygies. This method actually yields computational
methods, by using the geometry of the Cayley graph, and the
notion of deformation retraction of this universal cover.

Thus crossed complexes give a useful algebraic model of the
category of $CW$-complexes and cellular maps. This model does
lose a lot of information, but its corresponding advantage is
that it allows for algebraic description and computation, for
example of morphisms and homotopies. This is the key aspect of
the methods of \cite{BR}. See also the results in Theorems 3.4 -
3.6 here.

Many geometric and algebraic situations are specified by related,
and often more complicated,  non abelian information, not readily
computable by traditional means. As an example, we mention the
non abelian tensor product of groups due to Brown and Loday
\cite{B-Lo}.

Thus we can say that crossed complexes: \begin{enumerate}[(i)]
\item  give a first step to a full non abelian theory; \item  have
good categorical properties; \item  give a `linear' model of
homotopy types; \item  this model includes all homotopy 2-types;
\item are amenable to computation; \item  give one form of `higher
dimensional group'. \end{enumerate}

A further advantage of using crossed complexes of groupoids is
that this allows for the category $\crs$ to be monoidal closed:
there is a tensor product $-\otimes-$ and internal hom
$\Crs(-,-)$ such that there is a natural isomorphism
$$\crs(C \otimes D,E) \cong \crs(C,\Crs(D,E))$$
for all crossed complexes $C,D,E$. Here $\Crs(D,E)_0= \crs(D,E)$,
the set of morphisms $D \to E$, while $\Crs(D,E)_1$ is the set of
`1-fold left homotopies' $D \to E$. Note that while the tensor
product can be defined directly in terms of generators and
relations, such a definition makes it not easy to verify
essential properties of the tensor product, such as that the
tensor product of free crossed complexes is free. The proof of
this fact in \cite{BHclass} uses the above adjointness as a
necessary step to prove that  $-\otimes D$ preserves colimits.

An important result is that if $X_*,Y_*$ are filtered spaces,
then there is a natural transformation
$$\eta : \pi(X_*) \otimes \pi(Y_*) \to \pi(X_*\otimes Y_*)$$ which
is an isomorphism if $X_*,Y_*$ are $CW$-complexes (and in fact
more generally \cite{Ba-B}). In particular, the basic rules for
the tensor product are modelled on the geometry of the product of
cells $E^m \otimes E^n$ where $E^0$ is the singleton space, $E^1$
is the  interval $[-1,1]$ with two 0-cells and one 1-cell, while
for $m \ge 2$ $E^m=e^0 \cup e^{m-1} \cup e^m$. This leads to
defining relations for the tensor product. To give these we first
define   a bimorphism of crossed complexes (using additive
notation throughout).

A \emph{bimorphism} $\theta:(A,B) \to C$ of crossed complexes is
a family of maps $\theta:A_m \times B_n \to C_{m+n}$ satisfying
the following conditions, where $a \in A_m, b \in B_n$:

\begin{enumerate}
\item[\rm(i)] $\beta(\theta((a,b))=\theta(\beta a,\beta b)) \text{ for all }
a \in A, b \in B$.
\item[\rm(ii)] $\theta(a,b^{b_1})=\theta(a,b)^{\theta(\beta a,b_1)} \text{ if }
m \ge 0, n \ge 2$.
\item[\rm(ii)$'$] $\theta(a^{a_1},b)=\theta(a,b)^{\theta(a_1,\beta b)}\text{ if }
m \ge 2, n \ge 0$.
\item[\rm(iii)] \begin{align*}
\theta(a,b+b') & = \begin{cases} \theta(a,b) + \theta(a,b') &
\text{ if } m=0,n \ge 1 \  \text{or } m \ge 1,
n \ge 2, \\
\theta(a,b)^{\theta(\beta a,b')} + \theta(a,b') & \text{ if } m
\ge 1, n=1.
\end{cases}
\end{align*}
\item[\rm(iii)$'$] \begin{align*}
\theta(a+a',b) & = \begin{cases} \theta(a,b) + \theta(a',b) &
\text{ if } m \ge 1,n=0 \ \text{or } m \ge 2,
n \ge 1, \\
\theta(a',b) + \theta(a,b)^{\theta(a',\beta b)} & \text{if }
m=1,n \ge 1.
\end{cases}
\end{align*}
\item[\rm(iv)] \begin{align*}
\delta(\theta(a,b)) & = \begin{cases} \theta(\delta a,b) + (-)^m
\theta(a,\delta b) & \text{if } m \ge 2, n \ge 2,
\\
-\theta(a,\delta b) - \theta(\delta^1a,b) +
\theta(\delta^0a,b)^{\theta(a,
\beta b)} & \text{if } m=1,n \ge 2, \\
(-)^{m+1} \theta(a,\delta^1b) + (-)^m \theta(a,
\delta^0b)^{\theta(\beta a,
b)} + \theta(\delta a,b) & \text{ if } m \ge 2, n=1, \\
-\theta(\delta^1a,b) - \theta(a,\delta^0b) + \theta(\delta^0a,b)
+ \theta(a, \delta^1b) & \text{if } m=n=1.
\end{cases}
\end{align*}
\item[\rm(v)] \begin{align*}
\delta(\theta(a,b)) & = \begin{cases}
\theta(a,\delta b) & \text{if } m=0, n \ge 2, \\
\theta(\delta a,b) & \text{if } m \ge 2, n=0.
\end{cases}
\end{align*}
\begin{align*}
\delta^\alpha(\theta(a,b)) & = \begin{cases}
\theta(a,\delta^\alpha b) & \text{if }m=0,n=1 \ (\alpha=0,1), \\
\theta(\delta^\alpha a,b) & \text{if } m=1,n=0 \ (\alpha=0,1).
\end{cases}
\end{align*}
\end{enumerate}
The tensor product  of crossed complexes $A,B$ is then given by
the universal bimorphism $(A,B) \to A \otimes B$, $(a,b) \mapsto
a \otimes b$. So the rules for the tensor product are obtained by
replacing  $\theta(a,b)$ by $a \otimes b$ in the above.

\begin{example}
Let $\lan X|R \ran, \lan Y | S\ran$ be presentations of groups
$G,H$ respectively and let $C(R) \to F(X)$, $C(S) \to F(Y)$ be
the corresponding free crossed modules, regarded as crossed
complexes of length 2. Their tensor product $T$ is of length 4
and is given as follows:
\begin{itemize}
\item $T_1$ is the free group on generating set $X \sqcup
Y$;
\item $T_2$ is the free crossed $T_1$-module on $R \sqcup S \sqcup
(X \otimes Y)$ with the boundaries on $R,S$ as given before but
$$\delta_2(x \otimes y) = y\io x\io yx;$$
  \item $T_3$ is the free $(G \times H)$-module on generators $r\otimes y,
x \otimes s $ with boundaries
$$\delta_3(r \otimes y)= r \io r^y (\delta_2 r \otimes y ), \qquad
\delta_3(x \otimes s) = (x \otimes \delta_2s) \io s\io s^x;$$
  \item $T_4$ is the free $(G \times H)$-module on generators
$r\otimes s$, with boundaries
$$\delta _4(r \otimes s)= (\delta_2r \otimes s)+ (r \otimes \delta_2 s),$$
\end{itemize}
for $x \in X, y \in Y, r \in R, s \in S$.
\end{example}

The conventions here may seem (even are) awkward. They arise from
the derivation of the tensor product via another category of
`cubical $\omega$-groupoids with connections', and the formulae
are forced by our conventions for the equivalence of the two
categories \cite{BHalg,BHtens}. The important point is that we can
if necessary calculate with the formulae (i)-(v), because elements
such as $\delta_2 r \otimes y $ may be expanded using the rules
for the tensor product. Alternatively, the form $\delta_2 r
\otimes y $ may be left as it is since it naturally represents,
for example if $\dim y =1$, a subdivided cylinder.

A related result is that if $C,D$ are free crossed resolutions of
groups $C,D$ then $C \otimes D$ is a free crossed resolution of $G
\times H$, as proved by Tonks in \cite{Andy} and recovered in our
section 4. This allows for presentations of modules of identities
among relations for a product of groups to be read off from the
presentations of the individual modules. There is a lot of work on
generators for modules of identities (see for example \cite{HMS})
but not so much on higher syzygies.

As said above, the results of Brown and Higgins are not proved
directly, but use a category equivalent to $\crs$, namely a
category of `cubical $\omega$-groupoids with connections'. It is
in the latter category that the exponential law is easy to
formulate and prove, as is the construction of the natural
transformation $\eta$. However the proof of all the properties of
the equivalence is a long story.

In particular if we set $\I= \pi(\E^1)$, then a `1-fold left
homotopy' of morphisms $D \to E$ is defined to be a morphism $\I
\otimes D\to E$. The existence of this `cylinder object' $\I
\otimes D$ allows a lot of abstract homotopy theory \cite{Ka-P} to
be applied immediately to the category $\crs$. This is useful in
constructing homotopy equivalences of crossed complexes, using
for example gluing lemmas.

An important construction is the simplicial nerve $NC$ of a
crossed complex $C$. This is the simplicial set defined by
$$(NC)_n= \crs(\pi\Delta^n,C).$$ It directly generalises the nerve
of a group. In particular this can be applied to the internal hom
functor $\Crs(D,E)$ to give a simplicial set $N(\Crs(D,E))$ and
so turn the category $\crs$ into a simplicially enriched
category. This allows the full force of the methods of homotopy
coherence to be used \cite{CP}.

The {\it classifying space} $BC$ of a crossed complex $C$ is
simply the geometric realisation $|NC|$ of the nerve of $C$. This
construction generalises at the same time: the classifying space
of a group; an Eilenberg-Mac~Lane space $K(A,n),\ n \ge 2$; the
classifying space for local coefficients. It also includes the
notion of classifying space $B\cM$ of a crossed module $\cM=(\mu:
M \to P)$. Every connected $CW$-space has the homotopy 2-type of
such a space, and so crossed modules classify all connected
homotopy 2-types. This is one way in which crossed modules are
naturally seen as 2-dimensional analogues of groups.

\section{A Generalised Van Kampen Theorem}
This theorem states roughly that the functor $\pi : \ftop \to
\crs$ preserves certain colimits. This allows the calculation of
certain homotopically defined crossed complexes, and in particular
to see how free crossed complexes arise from $CW$-complexes.

The following definition with one modification is taken from
\cite{Bhcolimits}. The modification which we emphasise is that in
\cite{Bhcolimits} the assumption is made throughout that all
filtered spaces satisfy the condition, there called $J_0$, that
each loop in $X_0$ is contractible in $X_1$. It seems desirable to
bring this condition into the following connectivity condition.

 \label{gvktcrs}
\Env{defn}{{\em   A filtered space $X\sast$ is called {\em
connected} if  the following conditions $\phi (X, m)$ hold for
each $m \ge 0:$
 \renewcommand{\labelenumi}{}
\begin{enumerate}
\item $J_0$: each loop in $X_0$ is contractible in $X_1$; \item
$\phi (X, 0)$: If $j > 0,$ the map $\pi_0 X_0 \rightarrow \pi_0
X_j,$ induced by inclusion, is surjective \item $\phi (X, m), \;(m
\ges 1):$ If $j > m$ and $\nu \in X_0,$ then the map
$$\pi_m (X_m , X_{m-1} , \nu) \rightarrow \pi_m (X_j , X_{m-1} ,
\nu ) $$  induced by inclusion, is surjective.
 \end{enumerate} }}
The following result gives another useful formulation of this
condition. We omit the proof. \Env{prop}{ A filtered space $X$ is
connected if and only if it is $J_0$ and for all $ n > 0$ the
induced map $\pi_0 X_0 \rightarrow \pi_0 X_n$ is surjective and
for all $r
> n
> 0$ and $\nu \in X_0 , \pi_n (X_r , X_n , \nu ) = 0.$ }

The filtration of a $CW$-complex by skeleta is a standard example
of a connected filtered space.

Suppose for the rest of this section that $X\sast$ is a filtered
space. Let $X=X_{\infty}$.

We suppose given a cover ${\mathcal U} = \{ U^\lambda \}_{\lambda
\in \Lambda}$ of $X$ such that the interiors of the sets of
${\cal{U}}$ cover $X.$ For each $\zeta \in \Lambda^n$ we set
$$U^{\zeta} = U^{\zeta_{1}} \cap \cdots \cap U^{\zeta_{n}} ,
U^\zeta_i = U^\zeta \cap X_i.$$
 Then $U^{\zeta}_0 \subseteq U^{\zeta}_1\subseteq \cdots$ is called the
{\it induced filtration} $U^{\zeta}\sast$ of $U^{\zeta}$.
Consider the following $\pi$-{\it diagram} of the cover:
\begin{equation}{ \diagram  \bigsqcup_{\zeta \in \Lambda^{2}} \pi
U^{\zeta}\sast\; \rto<0.5ex>^a \rto<-0.5ex>_b &\;
\bigsqcup_{\lambda \in \Lambda} \pi  U^{\lambda }\sast\;
 \rto^(0.6)c & \;\pi X\sast \enddiagram}\end{equation}

\noindent Here $\bigsqcup$ denotes disjoint union (which is the
same as coproduct in the category of crossed complexes); $a, b$
are determined by the inclusions $a_\zeta : U^{\lambda} \cap
U^{\mu} \rightarrow U^{\lambda}, b_{\zeta} : U^{\lambda} \cap
U^{\mu} \rightarrow U^{\mu}$ for each $\zeta = (\lambda , \mu )
\in \Lambda^2$; and $c$ is determined by the inclusions
$c_{\lambda} : U^{\lambda} \rightarrow X.$

The following result constitutes a generalisation of the Van
Kampen Theorem for the fundamental group (or groupoid).
\begin{thm}{\rm   (The
coequaliser  theorem for crossed complexes: Brown and Higgins
\cite{Bhcolimits})}
 \label{coeqcxth} Suppose that for
 every finite intersection
$U^{\zeta}$ of elements of ${\mathcal{U}}$ the induced filtration
$U^{\zeta}\sast$ is connected.  Then
\begin{enumerate}\item[{\em (C)}] $X \sast$ is connected, and
\item[{\em (I)}]
 in the above $\pi$-diagram
of the cover, $c$ is the coequaliser of $a, b$ in the category of
crossed complexes. \end{enumerate} \label{crscoeq} \end{thm}

The proof of this theorem is not at all straightforward, and uses
another category equivalent to that of crossed complexes, called
the category of cubical $\omega$-groupoids with connections
\cite{BHalg}.  It is this category which is adequate for two key
elements of the proof, the notion of `algebraic inverse to
subdivision', and the `multiple compositions of homotopy addition
lemmas' \cite{Bhcolimits}. The setting up of this machinery takes
considerable effort.

In this paper we shall take as a corollary that the coequaliser
theorem applies to the case when $X$ is a $CW$-complex with
skeletal filtration and the $U^\lambda$ form a family of
subcomplexes which cover $X$.

In order to apply this result to free crossed resolutions, we
need to replace free crossed resolutions by $CW$-complexes. A
fundamental result for this is the following, which goes back to
Whitehead \cite{W-SHT} and Wall \cite{Wa}, and which is discussed
further by Baues in \cite[Chapter VI, \S 7]{Ba}:
\begin{thm}
Let $X_*$ be a $CW$-filtered space, and let $\phi: \pi X_* \to C$
be a homotopy equivalence to a free crossed complex with a
preferred free basis. Then there is a $CW$-filtered space $Y_*$
and an isomorphism $\pi Y_* \cong C$ of crossed complexes with
preferred basis such that $\phi$ is realised by a homotopy
equivalence $X_* \to Y_*$.
\end{thm}
In fact, as pointed out by Baues,  Wall states his result in terms
of chain complexes, but the crossed complex formulation seems more
natural, and avoids questions of realisability in dimension $2$,
which are unsolved for chain complexes.

\begin{cor}\label{cwmodel}
If $C$ is a free crossed resolution of a group $G$, then $C$ is
realised as free crossed complex with preferred basis by some
$CW$-filtered space $Y_*$.
\end{cor}
\begin{proof}
We only have to note that the group $G$ has a classifying
$CW$-space $BG$ whose fundamental crossed complex $\pi BG$ is
homotopy equivalent to $C$.
\end{proof}

Baues also points out in \cite[p.657]{Ba} an extension of these
results which we can apply to the realisation of morphisms of
free crossed resolutions.

\begin{prop}\label{cwmaps}
Let $X=K(G,1), Y=K(H,1)$ be $CW$-models of Eilenberg -Mac Lane
spaces and let $h:\pi(X_*) \to \pi(Y_*)$ be a morphism of their
fundamental crossed complexes with the preferred bases given by
skeletal filtrations. Then $h=\pi(g)$ for some cellular $g: X \to
Y$.
\end{prop}
\begin{proof}
Certainly $h$ is homotopic to  $\pi(f)$ for some $f:X \to Y$ since
the set of pointed homotopy classes $X \to Y$ is bijective with
the morphisms of groups $G \to H$. The result follows from
\cite[p.657,**]{Ba} (`if $f$ is $\pi$-realisable, then each
element in the homotopy class of $f$ is $\pi$-realisable').
\rule{0mm}{1mm}
\end{proof}

Note that from the computational point of view we will start with
a morphism $G \to H$ of groups and then lift that to a morphism of
free crossed resolutions. It is important for our methods that
such a morphism is exactly realised by a cellular map of the
cellular models of these resolutions. This is useful also because
a cellular map may be difficult to describe in geometric terms,
whereas it is not so hard to write down algebraically the lift of
a morphism $\sC_p \to \sC_{rp}$ to a morphism $F(\sC_p) \to
F(\sC_{rp})$, as is done in \cite{M}.

These results give a strategy of weaving between spaces and
crossed complexes. The key problem is to prove that a construction
on free crossed resolutions yields an aspherical free crossed
complex, and so also a resolution. The previous result allows us
to replace the free crossed resolutions by $CW$-complexes. We can
also replace morphisms of free crossed resolutions by cellular
maps.  We have a result of Whitehead \cite{W-asph} which allows us
to build up $K(G,1)$s as pushouts of other $K(G,1)$s provided the
induced morphisms of fundamental groups are injective. The
Coequaliser Theorem now gives that the resulting fundamental
crossed complex is exactly the one we want. More precise details
are given in the last section.

Note also an important feature of this method: {\it we use
colimits rather than exact sequences}. This enables precise
results in situations where exact sequences might be inadequate,
since they often give information only up to extension.

The relation of crossed complex methods to the more usual chain
complexes with operators is studied in \cite{BHchain}, developing
work of Whitehead \cite{W1}.

\section{Application to graph products of groups} Graph products of groups have
been studied for example in \cite{D,BHP,Co}. The paper \cite{BHP}
obtains generators for the module of identities among relations,
while \cite{Co} obtains a projective resolution of the graph
product given projective resolutions of the individual groups. Of
course our aim is for the information to be given in the
non-abelian form, where appropriate.

Let $\Gamma$ be a finite, undirected graph without loops or
multiple edges. Suppose the vertices of $\Gamma$ are well ordered,
and suppose given for each vertex $p$ of $\Gamma$ a group $G_p$.
The {\it graph  product} of the groups
$$G_\Gamma = \prod_\Gamma^p G_p$$ is obtained from the free product of the
groups $G_p$ (in the given order) by adding the relations
$[G_p,G_q]=\{1\}$ whenever $(p,q)$ is an edge of $\Gamma$.

We now consider the problem of constructing a free crossed
resolution of the graph product $G_\Gamma$ given free crossed
resolutions $C_p$ of each group $G_p$. To this end we define graph
products in two other contexts.

First of all we define the nerve $N(\Gamma)$ to be the simplicial
complex whose (ordered) simplices are the complete subgraphs of
$\Gamma$.

Next suppose we are given for each vertex $p$ of $\Gamma$ a
pointed crossed complex  $C_p$. For each simplex $\sigma=(p_0,
\ldots, p_n)$ of $N(\Gamma)$ let $C_\sigma$ be the subcomplex of
the tensor product of all the $C_p$ (in the given vertex order)
such that the $p$th coordinate is the base point of $C_p$ if $p
\not\in \sigma$ and is otherwise arbitrary. Let $C_\Gamma $ be the
subcomplex generated by  all the $C_\sigma$ for all $\sigma \in
N(\Gamma)$. This we call the {\it graph tensor product}
$$\otimes_\Gamma^p C_p$$ of the crossed complexes $C_p$.

Next suppose we are given for each vertex $p$ of $\Gamma$ a
pointed $CW$-space $X_p$. For each simplex $\sigma=(p_0, \ldots,
p_n)$ of $N(\Gamma)$ let $X_\sigma$ be the subset of the product
of all the $X_p$ (in the given vertex order) such that the $p$th
coordinate is the base point of $X_p$ if $p \not\in \sigma$ and
is otherwise arbitrary. Let $X_\Gamma $ be the union of all the
$X_\sigma$ for all $\sigma \in N(\Gamma)$. This we call the {\it
graph product} of the spaces $X_p$.

The main results on these spaces are: \begin{thm} If each space
$X_p$ is aspherical, then so also is their graph product.
\end{thm}
\begin{proof}
The proof is modelled on that given by Cohen in \cite{Co}.

If $\Gamma$ is a complete graph the result is clear, since finite
products of aspherical spaces are aspherical. We now work by
induction on the number of vertices of $\Gamma$.

Suppose $\Gamma$ is not complete.   Then there are vertices $p,q$
of $\Gamma$ which do not form an edge of $\Gamma$. Let $V$ be the
vertex set of $\Gamma$, and let $\Gamma_0,\Gamma_1,\Gamma_2$ be
the full subgraphs of $\Gamma$ on the complements in $V$ of
$\{p,q\}, \{p\}, \{q\}$ respectively. Let $X_0,X_1,X_2$ be the
corresponding graph products.  By the inductive assumption, these
are aspherical.  It  is clear that
$$X_\Gamma= X_1 \cup_{X_0} X_2.$$ But the maps on fundamental
groups induced by the inclusions $X_0 \to X_1, X_0 \to X_2$ are
injective. A theorem of Whitehead \cite{W-asph} now implies that
$X_\Gamma$ is aspherical.
\end{proof}

\begin{thm}
The fundamental crossed complex of the graph product of
$CW$-spaces is the graph tensor product of their fundamental
crossed complexes.
\end{thm}
\begin{proof}
This is immediate from previous results on the fundamental crossed
complex of a product of $CW$-spaces, and the Generalised Van
Kampen Theorem.
\end{proof}
\begin{cor}
If each $C_p$ is a free crossed resolution of a group $G_p$ then
the graph tensor product of the crossed resolutions is a free
crossed resolution of the graph product of the groups.
\end{cor}

Note that our result is stronger than that of Cohen in \cite{Co}
in that we obtain non abelian information. On the other hand he
obtains information on resolutions over arbitrary rings which
cannot currently be obtained by our methods.

\begin{example} In the case of a direct product $A\times B \times
C$ of groups, generators $x,y,z$ of $A,B,C$ respectively give rise
to an element of dimension $3$ of the corresponding tensor product
of free crossed resolutions, namely $x \otimes y \otimes z$. The
boundary of this is an identity among relations, and can be worked
out explicitly from the formulae for the tensor product. It in
fact corresponds to the cubical Homotopy Addition Lemma (HAL) in
\cite[Lemma 7.1]{BHalg}. This gives another view of Loday's
`favourite example' in \cite{Lo}, which is the case $A=B=C=\sC$,
the infinite cyclic group.
\end{example}
\begin{example}Suppose we take the graph product of four infinite cyclic groups
$A,B,C,D$ with generators $x,y,z,w$ and where the graph $\Gamma$
is a square with vertices assigned to the groups $A,B,C,D$ in that
order clockwise, say. Each group has a free crossed resolution of
length 1, with the same generators, say. The graph tensor product
of these resolutions then has free generators: $x,y,z,w$ in
dimension 1; $x \otimes y, y \otimes z, z \otimes w, x \otimes w$
in dimension 2; and no generators in dimension 3.
\end{example}

\end{document}